 

\def\P{{\cal P}}
\baselineskip=14pt
\parskip=10pt
\def\Tilde{\char126\relax}

\font\eightrm=cmr8  
\font\eighttt=cmtt8
\magnification=\magstephalf

\parindent=0pt
\overfullrule=0in
 
\bf
\centerline
{ENUMERATION SCHEMES AND, MORE IMPORTANTLY, THEIR AUTOMATIC GENERATION}
\bigskip
\centerline{ {\it Doron ZEILBERGER}\footnote{$^1$}
{\eightrm  \raggedright
Department of Mathematics, Temple University,
Philadelphia, PA 19122, USA. 
{\eighttt zeilberg@math.temple.edu} \hfill \break
{\eighttt http://www.math.temple.edu/\Tilde zeilberg   .}
May 13, 1998. Supported in part by the NSF.
} 
}
\rm 
{\bf Towards the Virtual Combinatorialist}
 
It is way too soon to teach our computers how to become
full-fledged {\it humans}. It is even premature to teach them
how to become {\it mathematicians}, it is even unwise, at present, to teach
them how to become {\it combinatorialists}. But the time
is ripe to teach them how to become  experts in a suitably
defined and narrowly focused subarea of combinatorics.
In this article, I will describe my efforts to teach
my beloved computer, Shalosh B. Ekhad, how to be an enumerator
of Wilf classes. 
 
Just like the Wright Brothers' first flight, and Len Adelman's
first DNA computer (and perhaps  like the first Fleischmann-Pons
Cold Fusion experiment), the actual accomplishments are very meager. But
I do hope to demonstrate {\it feasibility}, and 
{\it jump-start} a research area that is not quite part of AI, and not
quite algorithmic combinatorics. It is something
brand-new, let's call it {\it Artificial Combinatorics} (AC). 
The closest analogy is
Deep Blue, but instead of playing Chess, Shalosh will play
a game called `Wilf-class enumeration'. It is still very far
from beating the Kasparovs of the area 
(Miklos Bona, Alex Burstein, John Noonan, Frank Schmidt, Rodica Simion, 
Zved Stankova, Julian West, to name a few), but 
I am sure that the ultimate version will. 
 
With all due respect to Wilf classes and Enumeration, and even
to Combinatorics, the main point of this article is not to enhance
our understanding of Wilf classes, but to {\it illustrate} how
much (if not all) of mathematical research will be conducted 
in a few years. It goes as follow. Have a (as of now, human)
mathematician get a brilliant idea. Teach that idea to a computer,
and let the computer `do research' using that idea. 
 
{\bf Accompanying Software}
 
In a sense, this article is a user's manual for
the Maple package {\tt WILF}, available from
my homepage {\tt http://www.math.temple.edu/\Tilde zeilberg}
(click on {\tt Maple packages and programs}, and then download
{\tt WILF}). The empirical version of the rigorous
{\tt WILF} is {\tt HERB}, also available there.
 
{\bf Enumeration Schemes}
 
Suppose that we have to find a `formula' (in the sense of Wilf[Wi],
i.e. a polynomial-time algorithm) for computing $a_n:=\vert A_n\vert$,
where $A_n$ is an infinite family of finite sets, parameterized
by $n$. Usually $A_n$ is a natural subset of a larger set
$B_n$, and is defined as the set of members of $B_n$ that
satisfy a certain set of conditions $C_n$. For example
if $A_n$ is the set of permutations on $\{1,2, \dots , n\}$,
then $B_n$ may be taken as the set of words of length
$n$ over the alphabet $\{1,2, \dots , n\}$, and $C_n$ can
be taken as the condition: `no letter can appear twice'.
A naive algorithm for enumerating $A_n$ would be to
actually {\it construct} the set, by examining the
members of $B_n$, one by one, checking whether they
satisfy $C_n$, and admitting those that qualify.
Then $a_n$=Cardinality of $A_n$.
 
But a much better approach would be to find a
{\it structure} theorem that expresses $A_n$, using unions,
Cartesian products, and possibly complements, of well known sets.
Failing this, it would be also nice to express $A_n$,
recursively, in terms of $A_{n-1}, A_{n-2}, \dots $, and easy-to-count
sets, getting a {\it recurrence formula}. Going back to the
permutation example, Levi Ben Gerson ([L]) proved the structure
theorem $A_n \equiv \{1,2, \dots , n\} \times A_{n-1}$, from which
he deduced the {\it recurrence} $a_n=n a_{n-1}$, enabling a
polynomial-(in fact linear-) time algorithm for computing
$a_i$, for $1 \leq i \leq n$.
 
Alas, this is not always easy, and for many enumeration sequences,
e.g. the number of self-avoiding walks, may well be
{\it impossible}, and who knows, perhaps one day even
{\it provably impossible}.
 
It is conceivable, however, that a combinatorial family
$A(n)$, does not possess a recursive structure by itself, but by
{\it refining it}, using a suitable parameter, one can partition
$A(n)$ into the disjoint union:
$$
A(n)= \bigcup_{i=1}^{n} B(n,i) \quad, 
$$
and try and find a {\it structure theorem} for the two-parameter
family $B(n,i)$. This will imply a recurrence for the
cardinalities $b(n,i):=\vert B(n,i) \vert$, that would enable
a fast algorithm for $a(n)=\sum_{i=1}^{n} b(n,i)$.
 
Sometimes, not even the $B(n,i)$ suffice. Then we could try to
partition $B(n,i)$ into the following disjoint union:
$$
B(n,i)= \bigcup_{j=1}^{i-1} C_1(n,i,j) \cup \bigcup_{j=i+1}^{n} 
C_2(n,i,j)  \quad,
$$
and try to express $C_1(n,i,j)$ and $C_2(n,i,j)$ in terms of
$A(m)$, $B(m,i')$, $C_1(m,i',j')$, and $C_2(m,i',j')$, with
$m<n$. One can keep going {\it indefinitely}.
If this process {\it halts} after a finite number
of refinements, then we have indeed a {\it formula} (in the
sense of Wilf) for $a(n)$.
 
{\bf Permutations with Forbidden Patterns; Wilf classes}
 
A {\it pattern} of length $k$ is a permutation of
$\{ 1, 2, \dots, k \}$.
 
The {\it reduction } of a vector of distinct integers of length $k$
is the pattern obtained by replacing the smallest
element by $1$, the second smallest by $2$, $\dots$,
the largest by $k$. For example, the reduction of
$264$ is $132$.
 
A permutation $\pi$ of length $n$, of an ordered set, is said
to contain the pattern $\sigma$ of length $k$ if there
are {\it places} $i_1 < \dots < i_k$, such that
$\pi_{i_1} \pi_{i_2} \dots \pi_{i_k}$ reduces to $\sigma$.
For example $51872463$ contains the pattern $3421$ 
(with $i_1=1, i_2=3, i_3=6, i_4=8$).
 
Our Goal is to investigate the following
 
{\bf Problem:} Given a set of patterns $\P$, study the
sequence $a(n;\P):=\vert A(n;\P)\vert$,
where $A(n;\P)$ is the set of permutations on
$\{1,2, \dots, n \}$ that avoid the patterns of $\P$.
 
We will sometimes write $A(n)$ and $a(n)$ instead of
$A(n;\P)$ and $a(n;\P)$, where the $\P$ is implied
from the context.
 
Of course $a(n;\emptyset) \equiv n!$, and 
$a(n;\{12\})\equiv 1$. It is well known (see [We][K]) that
$a(n;\{123\})=a(n;\{132\})=C_n$, where $C_n=(2n)!/(n!(n+1)!$ are the
Catalan numbers.
 
Naively, it takes an exponential time to compute 
the first $n$ terms of the sequence $a(n;\P)$. The best
possible scenario is an explicit formula for the sequence.
Failing this, it would be almost as nice to have a linear recurrence
equation with polynomial coefficients (and hence inducting it
in the hall of fame of P-recursive sequences).
 
Sometimes it happens that there exist two sets of patterns
$\P$ and $\P'$ such that $a(n;\P)=a(n;\P')$ for all $n$.
If this happens for a non-trivial reason, we say that
$\P$ and $\P'$ belong to the same {\it Wilf class}.
 
The approach that I am going  to investigate here, and that
I taught Shalosh B. Ekhad, is to partition the
set $A(n; \P)$ into a finite disjoint union of subsets, indexed by
certain prefixes, and to try to deduce {\it rigorously},
certain recurrence relations between them, that would enable
one to set up an {\it enumeration scheme}, that would lead to
a polynomial-time algorithm for computing $a(n; \P)$.
In the future it is hoped that these schemes could also be used
to derive (once again rigorously) generating functions for the
more refined quantities, that should lead to a constant-term
expression for $a(n;\P)$, and that would make it decidable
whether $\P$ and $\P'$ belong to the same Wilf class.
 
{\bf Apology:} The success rate of the present method, in its
present state, is somewhat disappointing. Ekhad was able to reproduce
the classical cases, and a few new ones, but for most
patterns and sets of patterns, it failed to find a scheme
(defined below) of a reasonable depth. But the present
set-up for setting up a scheme could be modified and extended
in various ways. We do believe that an appropriate
enhancement of the present method would yield, if not
a one-hundred-per-cents success rate, at least close to it.
 
{\bf The Most Trivial Non-Trivial Example}
 
In order to illustrate the method, let's work out {\it in full
detail}, and in {\it plain English}, an Enumeration Scheme for
the set of permutations avoiding the pattern $123$.
 
Let $A(n)$ be the set of such permutations. Break $A(n)$ up into
the union
$$
A(n)=\cup_{i=1}^{n} A_1(n;i) \quad ,
$$
where $A_1(n,i)$ are the set of $123$-avoiding permutations 
on $[1,n]$ that start with $i$, i.e.
$$
A_1(n,i):=\{ \pi \in A(n) ; \pi_1=i \} \quad .
$$
 
Let's break up $A_1(n;i)$ into the following disjoint union:
$$
A_1(n,i)=\bigcup_{h=1}^{i-1} A_{21}(n;h,i) \cup 
\bigcup_{j=i+1}^{n} A_{12}(n;i,j) 
\quad,
$$
where for $1 \leq i < j \leq n$,
$$
A_{12}(n,i,j):=\{ \pi \in A(n) ; \pi_1=i \,,\, \pi_2=j \} \quad ,
$$
$$
A_{21}(n,i,j):=\{ \pi \in A(n) ; \pi_1=j\,,\, \pi_2=i \} \quad .
$$
 
Now let's look at a typical element $\pi$ of $A_{21}(n,i,j)$.
It has $\pi_1=j$ and $\pi_2=i$. Deleting the first entry
(i.e. the $j$), and reducing, gives us a member $\sigma$ of
$A_1(n-1,i)$, since you can't create {\it trouble} (i.e. a $123$ pattern)
by skipping town. FURTHERMORE, if you start out with {\it any}
element $\sigma$ of $A_1(n-1,i)$, make room for the $j$ (by raising all
members $\geq j$ by $1$), and insert $j$ at the very beginning,
you can't create new trouble, because the $j$ is
``body-guarded'' by the $i$. Let me explain.
If inserting $j$ at the very beginning
would have caused an undesirable $123$-pattern, then of course,
the ``$1$'' of the ``$123$'' would be $j$, but the ``$2$'' can't be $i$ (since
$j>i$), hence in this hypothetical new $123$ that the insertion of
$j$ at the very beginning would have created,
the ``$2$'' would have been a certain $\sigma_a$ and the ``$3$'' 
would be a certain $\sigma_b$, with $2 <a<b$. But had that been the case, then
$i$ would have been an EVEN BETTER ``$1$''!, because $i<j$,
and $\sigma_1<\sigma_a<\sigma_b$ would have already been a $123$-pattern.
It follows that $A_{21}(n;i,j)$ is isomorphic to $A_1(n-1;i)$.
 
Now let's examine $A_{12}(n;i,j)$. If $j<n$, then $n$ must be
somewhere to the right, making $ijn$ a $123$-pattern. Hence
$A_{12}(n;i,j)$ is empty when $j<n$. On the other hand, if
$j=n$, then deleting it creates a member of $A_1(n-1;i)$.
Conversely, given any member of $A_1(n-1;i)$, inserting
$n$ right after the first entry (hence making $n$ the
second entry), can't create a $123$. Hence 
$A_{12}(n;i,n) \equiv A_1(n-1;i)$. Summarizing, we have
the following structure:
$$
A(n)=\bigcup_{i=1}^{n} A_1(n;i) \quad ,
$$
$$
A_1(n,i):=\bigcup_{h=1}^{i-1} A_{21}(n;h,i) \cup 
\bigcup_{j=i+1}^{n} A_{12}(n;i,j) 
\quad,
$$
$$
A_{21}(n;i,j) \equiv A_1(n-1;i) \quad ,
$$
$$
A_{12}(n;i,j) \equiv \cases{\emptyset,& if $j < n$;\cr
             A_1(n-1;i),& if $j=n$.\cr}
$$
 
Now let $a(n):=\vert A(n) \vert$,  $a_1(n;i):=\vert A_1(n;i) \vert$,
$a_{12}(n;i,j):=\vert A_{12}(n;i,j) \vert$,
$a_{21}(n;i,j):=\vert A_{21}(n;i,j) \vert$.
Taking cardinalities gives us the following scheme:
$$
a(n)=\sum_{i=1}^{n} a_1(n;i) \quad ,
$$
$$
a_1(n,i):=\sum_{h=1}^{i-1} a_{21}(n;h,i) + \sum_{j=i+1}^{n} a_{12}(n;i,j) 
\quad,
$$
$$
a_{21}(n;i,j) = a_1(n-1;i)
$$
$$
a_{12}(n;i,j) =\cases{0,& if $j < n$;\cr
             a_1(n-1;i),& if $j=n$.\cr}
$$
Since $a_{21}$ and $a_{12}$ can be expressed in terms of $a_1$, the above
can be condensed to 
$$
a(n)=\sum_{i=1}^{n} a_1(n-1;i) \quad ,
$$
$$
a_1(n,i)= \left ( \sum_{h=1}^{i-1} a_1(n-1;h) \right ) +  a_1(n-1;i)=
\sum_{h=1}^{i} a_1(n-1 \,;\, h) \quad,
$$
from which it follows that $a_1(n;i)-a_1(n;i-1)=a_1(n-1;i)$.
This recurrence and the obvious boundary conditions,
$a_1(n;1)=1$ and $a_1(n;n+1)=0$, yield the explicit solution
$a_1(n;i)={{n+i-2} \choose {n-1}} -{{n+i-2} \choose {n}}$
Hence $a(n)=a_1(n+1,n+1)=C_n={{2n} \choose {n}}/(n+1)$.
 
{\bf Introducing Prefix Schemes}
 
The above example leads to the following definitions.
Fix a set of patterns $\P$.
Let $\sigma=\sigma_1 \dots  \sigma_k$ be any permutation
of length $k$. For $1 \leq i_1 < i_2 < \dots < i_k \leq n$,
let
$$
A_\sigma(n; \P; i_1, \dots , i_k) :=
\{ \pi \in A(n; \P) \,\, \vert \,\, \pi_1=i_{\sigma_1}, \pi_2=i_{\sigma_2},
\,\,\dots \,\,, \,\,\pi_k=i_{\sigma_k} \} \quad .
$$
For example $A_{132}(5; \{1234,1432\}; 2,3,5) := \{25314, 25341 \}$.
 
The set of {\it refinements} of a permutation 
$\sigma=\sigma_1 \dots \sigma_k$ of lengths $k$ is the set of $k+1$ 
permutations of length $k+1$ that have the property that deleting their last 
entry reduces to $\sigma$. We will denote the set of refinements of $\sigma$ by
$Refinements(\sigma)$. For example, $Refinements(312)=\{ 3124,4123,4132,4231\}$.
 
For any set of patterns $\P$, and any prefix-permutation 
$\sigma$, we have the obvious recurrence
$$
A_\sigma(n;\P;i_1, \dots , i_k)=
\bigcup_{r=1}^{k+1} \bigcup_{j=i_{r-1}+1}^{i_r-1}
A_{\sigma^{(r)}}(n;\P;i_1, \dots , i_{r-1}, j, i_r, \dots ,i_k) \quad ,
\eqno(*)
$$
where $\sigma^{(r)}$ is that element of $Refinements(\sigma)$ that ends with
$r$. We agree that $i_0=0$ and $i_{k+1}=n$. For example
$$
A_{132}(n;i_1,i_2,i_3)=
\bigcup_{j=1}^{i_1-1} A_{2431}(n; j,i_1,i_2,i_3) \,\,\cup
\bigcup_{j=i_1+1}^{i_2-1} A_{1432}(n; i_1,j,i_2,i_3) \,\,\cup
\bigcup_{j=i_2+1}^{i_3-1} A_{1423}(n; i_1,i_2,j,i_3) \,\, \cup
$$
$$
\bigcup_{j=i_3+1}^{n} A_{1324}(n; \P; i_1,i_2,i_3,j) \,\,  \quad .
$$
 
Given a set of patterns, $\P$, and a prefix permutation $\sigma$,
it may happen that all the members of $A_\sigma(n;\P)$, for all $n$,
must have either one or more of the properties $i_1=1$, or $i_1=i_2$,
$i_2=i_3, \dots, i_k=n$. Let $J(\sigma)$ be the set of
$0\leq r \leq k$ such that $i_r=i_{r+1}$ is {\it forced} for all
members of $A_\sigma(n;\P)$. In the above example with $\P=\{123\}$,
we had $J(12)=\{2\}$, since every $123$-avoiding permutation
that starts with $i_1i_2$, where $i_1<i_2$, {\it must} have $i_2=n$.
{\bf Another Example:} If $\P=\{1234,1324,1243\}$,
then $J(2413)=\{4\}$, since every member of
$A_{2413}(n;\P;i_1,i_2,i_3,i_4)$ must have its second entry, $i_4$ be
equal to $n$, or else the $n$ would be to the right of the fourth entry,
that would cause $i_2 i_4 i_3 n$ to be an illegal $1324$.
 
Consider all permutations. Let's remember our goal.
Given a set of patterns $\P$, we want to count the number
of permutations that avoid all the patterns of $\P$.
Let's call these permutations law-abiding. Hence a permutation
that has one or more occurrences of the patterns $\P$ should
be called {\it criminal}. Obviously, deleting any entry of
a  law-abiding permutation (and reducing), gives rise to
another law-abiding permutation. Hence given a prefix
permutation $\sigma=\sigma_1 \dots \sigma_k$, then for each place
$r$ ($1 \leq r \leq k$), we have the obvious inclusion:
$$
F_r(A_\sigma(n;\P;i_1, \dots , i_k)) \subset
A_{\sigma_{(r)}}(n-1;\P;i_1, \dots, i_{r-1},i_{r+1}-1, \dots , i_k-1) \quad ,
\eqno(Rodica)
$$
where $F_r$ is the injective mapping that consists of deleting
$i_r$ (wherever it is), and reducing, and $\sigma_{(r)}$ is the
permutation obtained from $\sigma$ by deleting $r$ and reducing.
 
BUT, it may happen that the inclusion $(Rodica)$ is an {\it equality}.
This leads to the following
 
{\bf IMPORTANT DEFINITION}. Given a set of patterns $\P$, and
a prefix permutation $\sigma$, the place $r$ ($1 \leq r \leq k$) is
{\it reversely deleteable} if for all $1 \leq i_1 < \dots < i_k \leq n$,
$$
F_r(A_\sigma(n;\P;i_1, \dots , i_k))=
A_{\sigma_{(r)}}(n-1;\P;i_1, \dots, i_{r-1},i_{r+1}-1, \dots , i_k-1) \quad .
\eqno(Julian)
$$
 
We must find a way to know,
{\it a priori}, by using logical reasoning,
whether a place is reversely deleteable.
 
How can we be sure that inserting $i_r$ in any permutation of
$A_{\sigma_{(r)}}(n-1;\P;i_1, \dots, i_{r-1},i_{r+1}-1, \dots , i_k-1)$,
in the appropriate place, is a safe thing to do? We must be
assured that its insertion cannot cause any trouble. How
can we be sure? We have to look at all the conceivable
{\it events} that its insertion can cause, i.e. all the possible
forbidden patterns of $\P$, with one of its entries coinciding
with the inserted $i_r$. If for each such possible scenario,
in which $i_r$ participates,
we can {\it logically} deduce the existence of {\it another event},
in which $i_r$ {\it does not} participate,
then we are safe, since {\it we are covered}.
 
{\bf Example}: Let $\P=\{123\}$. For $\sigma=21$, the place $r=1$
is reversely deleteable, since suppose that inserting the
$i_2$ would have created a 123. Then there exist $a$ and $b$
in $\pi$ such that $i_2 a b$ is a $123$- pattern. Had that been the
case, then {\it kal vakhomer}, $i_1 a b$ is also a $123$-pattern.
 
{\bf Another Example}: Let $\P=\{1234,1324,1243\}$.
$\sigma=2413$. We saw above that $J(\sigma)=\{4\}$, i.e. the
permutations $\pi$ that we are examining  start with
$i_2 n i_1 i_3$ for some $1 \leq i_1<i_2<i_3 \leq n$.
I claim that the first place, i.e. the $i_2$ is reversely
deleteable. 
 
Let's look at the kind of trouble the insertion of $i_2$ at the
beginning could cause.
 
Event 1: There exist $4<b<c<d\leq n$ such that
$i_2 \pi_a \pi_b \pi_c$ is a delinquent $1234$.
Then $i_2$ can be bailed out by $i_1$, since
$i_1 \pi_a \pi_b \pi_c$ would have been a bad event
even before the insertion.
 
Event 2: There exist $4<b<c\leq n$ such that
$i_2 i_3 \pi_b \pi_c$ is a delinquent $1234$.
Then $i_2$ can be bailed out by $i_1$, since
$i_1 i_3 \pi_b \pi_c$ would have been a bad event
even before the insertion.
 
Event 3: There exist $4<a<b<c \leq n$ such that
$i_2 \pi_a \pi_b \pi_c$ is a delinquent $1324$.
Then $i_2$ can be bailed out by $i_1$, since
$i_1 \pi_1 \pi_b \pi_c$ would have been a bad event
even before the insertion.
 
Event 4: There exist $4<b<c \leq n$ such that
$i_2 i_3 \pi_b \pi_c$ is a delinquent $1324$.
Then $i_2$ can be bailed out by $i_1$, since
$i_1 i_3 \pi_b \pi_c$ would have been a bad event
even before the insertion.
 
Event 5: There exist $4<a<b<c \leq n$ such that
$i_2 \pi_a \pi_b \pi_c$ is a delinquent $1243$.
Then $i_2$ can be bailed out by $i_1$, since
$i_1 \pi_1 \pi_b \pi_c$ would have been a bad event
even before the insertion.
 
Event 6: There exist $4<b<c \leq n$ such that
$i_2 i_3 \pi_b \pi_c$ is a delinquent $1243$.
Then $i_2$ can be bailed out by $i_1$, since
$i_1 i_3 \pi_b \pi_c$ would have been a bad event
even before the insertion.
 
So any law-abiding permutation $\pi$, that
starts with  a triplet that reduces to
$312$ has an {\it a priori guarantee} that sticking
an entry at the very front, 
that make it have prefix that reduces to $2413$, will
not get the permutation $\pi$ in trouble.
Hence we can know for sure that, with $\P=\{1234,1324,1243\}$,
$$
A_{2413}(n; \P; i_1,i_2,i_3,i_4)=
\cases{\emptyset,& if $i_4 < n$;\cr
A_{312}(n-1 \,;\, \P; i_1,i_3-1,i_4-1),& if $i_4=n$.\cr} \quad.
$$
 
{\bf Two `Trivial' Examples}
 
Let's rederive Levi Ben Gerson's famous recurrence for
$a(n;\emptyset)$, i.e. $a(n;\emptyset)=na(n-1;\emptyset)$.
 
First, let's refine $A(n; \emptyset)$ into
$$
A(n;\emptyset)=\cup_{i_1=1}^{n} A_1(n;\emptyset;i_1) \quad .
$$
Now, note that $r=1$ in $\sigma=1$ is reversely deleteable, since
inserting an $i_1$ at the front can't cause any trouble,
since nothing is forbidden. Hence
$$
A_1(n;\emptyset;i_1) \equiv A(n-1;\emptyset) \quad,
$$
from which follows that
$$
a(n;\emptyset)=\sum_{i_1=1}^{n} a_1(n;\emptyset;i_1)=
\sum_{i_1=1}^{n} a(n-1;\emptyset)=n a(n-1;\emptyset) \quad .
$$
 
Next, let's find a scheme when $\P=\{12\}$, that is
let's try to compute $a(n;\{12\})$, the number of
permutations on $n$ objects with no occurrence of the
pattern $12$, i.e. the number of decreasing permutations.
First, as always
$$
A(n;\{12\})=\bigcup_{i_1=1}^{n} A_1(n;\{12\};i_1) \quad .
$$
Now $J(1)=\{2\}$, since any $12$-avoiding permutation must
have $i_1=n$, or else $i_1 n$ would be a forbidden pattern.
Also $r=1$ is a reversely deleteable place, since
inserting $i_1=n$ at the very beginning of a permutation in
$A(n-1;\{12\})$ can't possibly cause trouble. We hence have
the following recurrence:
$$
A(n;\{12\}) \equiv A(n-1;\{12\}) \quad,
$$
which implies that $a(n;\{12\})= a(n-1;\{12\})$, and since
$a(0;\{12\})=1$, we have proven, rigorously, that
$a(n;\{12\})=1$, for all $n>0$.
 
A {\it prefix scheme} can be defined independently of sets
of forbidden patterns $\P$ as follows.
 
{\bf The Formal Definition of Prefix Scheme}
 
{\bf Definition}: A {\bf PREFIX SCHEME} is a five-tuple
$[Redu,Expa,A,B,C]$ such that:
 
(i) $Redu$ is a finite set of permutations (of various lengths).
 
(ii) $Expa$ is another finite set of permutations (of various lengths).
It must include the empty permutation $\phi$. $Expa$ and $Redu$ are
disjoint.
 
(iii) $A$ is a table assigning to each member $\sigma=\sigma_1 \dots \sigma_k$ 
of $Redu$, a place $A[\sigma]$ in $\sigma$, i.e. an integer
$i$ between $1$ and $k$.
 
(iv) $B$ is a table assigning to each permutation 
$\sigma=\sigma_1 \dots \sigma_k$ of $Expa$, its $k+1$ refinements.
 
(v) $C$ is a table that assigns to each member 
$\sigma=\sigma_1 \dots \sigma_k$ of
$Redu\cup Expa$ a certain set $J(\sigma)$, which is a (possibly empty) 
subset of $\{0,1, 2, \dots , k \} $.
 
(vi) For each $\sigma \in Expa$, each member $B[\sigma]$ must belong
to $Expa \cup Redu$.
 
For example, the following $[Redu, Expa,A,B,C]$
is a prefix scheme: $Redu=\{12,21\}$, $Expa=\{\phi,1\}$,
$A[12]=1, A[21]=1$, $B[\phi]=\{1\}, B[1]=\{12,21\}$,
$C[12]=\{2\}, C[21]=C[1]=C[\phi]=\emptyset$.

We are now ready to interface the abstract notion of {\it prefix
scheme} to that of {\it forbidden patterns}.
 
{\bf Another Important Definition}: 
 
A {\bf Prefix Scheme} $[Redu, Expa, A,B,C]$ is an {\bf Enumeration
Scheme} for a set of forbidden patterns $\P$, if the following
conditions hold:
 
(i) For every $\sigma \in Redu$, let $\sigma'$ be the reduction of
the permutation obtained from $\sigma$ by deleting the entry at the
place specified by $A[\sigma]$. The following property holds:
If you take any permutation $\pi$ that avoids the patterns
in $\P$, and that has a prefix that reduces to $\sigma'$, and
insert a new entry right before the $(A[\sigma]+1)^{th}$ place
of $\pi$, in such a way that the new permutation has a prefix
that reduces to $\sigma$, then that new permutation also avoids
the patterns of $\P$. In other words, the insertion is always
safe.
 
(ii) Suppose that a permutation $\pi$ 
of $\{1,2, \dots, n \}$, avoids the patterns of $\P$,
and has a prefix that reduces to one of the permutations 
$\sigma=\sigma_1 \dots \sigma_k$ of $Expa \cup Redu$, and
the first $k$ entries of $\pi$ are 
$i_{\sigma_1} i_{\sigma_2} \dots i_{\sigma_k}$, where, of course,
$1 \leq i_1< i_2 < \dots < i_k \leq n$. Then
if $0 \in C[\sigma]$, we MUST have $i_1=1$, and if
$k \in C[\sigma]$, we MUST have $i_k=n$, and for any other
member $j$ of $C[\sigma]$, we MUST have $i_j=i_{j+1}$.
 
{\bf Implementing the Scheme}
 
Once the computer (or, in simple cases, the human) has found
a prefix-scheme for a set of patterns $\P$, then we have
a ``formula'' (in the Wilfian sense discussed above) for
enumerating $a(n;\P)$. 
 
If $\sigma=\sigma_1 \dots \sigma_k \in Expa$, 
let $\sigma^{(j)}$ ($1 \leq j \leq k+1$), be that member of
$B[\sigma]$ that ends with $j$. We have, for $1 \leq i_1 < \dots < i_k \leq n$
$$
a_\sigma(n;\P;i_1, \dots , i_k)=
\sum_{j=1}^{k+1} 
\sum_{r=i_{j-1}+1}^{i_j-1}
a_{\sigma^{(j)}}(n;\P;i_1, \dots, i_{j-1}, r, i_j,\dots , i_k) \quad,
$$
where $i_0=0$ and $i_{k+1}=n+1$.
 
If $\sigma=\sigma_1 \dots \sigma_k$ is in $Redu$, 
let $A[\sigma]=r$. 
If $C[\sigma]=\emptyset$, then
$$
a_\sigma(n;\P;i_1, \dots , i_k)=
a_{\sigma_{(r)}}(n-1;\P;i_1, \dots , i_{r-1}, i_{r+1}-1, \dots , i_k-1)
\quad ,
\eqno(Miklos)
$$
where $\sigma_{(r)}$ is the permutation obtained 
by reducing the permutation obtained from
$\sigma$ by deleting $\sigma_r$.
 
If $J:=C[\sigma] \neq \emptyset$ then $(Miklos)$ holds
whenever $(i_1, \dots, i_k)$ ``obeys $J$'', i.e.
whenever $i_j=i_{j+1}$ for all $j\in J$
(recall that $i_0=1$, $i_{k+1}=n$).
Otherwise $a_\sigma(n;\P;i_1, \dots , i_k)=\emptyset$.
 
Since the sets $Redu$ and $Expa$ are finite, the enumeration
scheme is well-defined.
 
{\bf How to Find a Scheme}
 
First a warning: there is no guarantee that a set of
patterns $\P$ possesses a finite scheme. In fact, based
on empirical evidence, most don't (and if they did, their
depth would be so large that the `polynomial' in the
``polynomial growth guarantee'' would be of such a high degree
as to make it almost as bad as exponential growth.
 
But that's what so nice about non-trivial human research. 
If we know {\it beforehand} that we are guaranteed to succeed, than
it is not research, but doing chores. 
So our ``algorithm'' is not known to halt. To make it a genuine
algorithm, we make the maximal depth part of the input, and
restrict the search to Prefix Schemes of bounded depth.
If we fail, then the program returns $0$.
 
The algorithm for looking for a Prefix Scheme for
a given set of patterns $\P$ is a formalization of the
procedure described above. The input is
a set of patterns, and an integer, MaximalDepth.
 
We start with the empty permutation $\phi$ as belonging to
$Expa$. Its only refinement is $1$. Unless
$\P=\emptyset$, $1$ would also belong to $Expa$. Its
refinements are $12$ and $21$. We examine each prefix 
permutation $\sigma$, in turn, and see whether it has
a non-empty $J:=C[\sigma]$ (forced relations), and equipped
with that $J$, we look for a reversely deleteable place.
This we do by looking at all the conceivable events
that the place under consideration can participate in,
and look at all the {\it implied events}, hoping to
see amongst them an event that does not include the
examined place. If indeed  each and every possible calamity
in which the examined place participates in, implies
another event in which it does not participate in,
(by using the transitivity of the order relation), then
that place is indeed reversely deleteable.
 
If such a reversely deleteable
place exists, then $\sigma$ becomes a member of $Redu$,
and the lucky place, that made it a member is recorded as
$A[\sigma]$. If there are no such places, we reluctantly
put $\sigma$ in $Expa$, and find its refinements, that
we store as $B[\sigma]$. We then examine these in turn, until
we either get a prefix-permutation of length $MaximalDepth+1$,
in  which case we sadly exit with $0$, or all the offsprings of
the member of $Expa$ are in $Expa \cup Redu$.
 
{\bf Important Remark:} All the above deductions are made
completely automatically by the computer. So we have
(a very primitive) `Artificial Intelligence'.
 
{\bf Using the Maple package} {\tt WILF}
 
The procedure in the package {\tt WILF}, that looks for 
Schemes for a set of patterns, is {\tt Scheme}.
The syntax is: ``Scheme(Set\_Of\_Patterns,Maximal\_Depth)''.
For example, to get a scheme for $\P=\{123,132\}$, type:
``{\tt Scheme($\{$[1,2,3],[1,3,2]$\}$,2);}'' (without the quotes, of course),
followed, of course, by Carriage Return. Having found
the scheme (let's call it {\tt sch}), to find
the number of permutations on $n$ objects avoiding the
patterns $\P$, you do ``$Miklos(n,sch);<CR>$''.
For example ``{\tt Miklos(3,sch);<CR>}'' should yield
{\tt 4}. To get the 
first $L$ terms of the sequence $a(n; \P)$ (after
$sch:=Scheme(\P;Depth)$ was successful for some $Depth$),
do ``{\tt SchemeSequence(sch,L);<CR>}''. The
package {\tt WILF} could also try to guess (empirically,
but perhaps rigorizably), a linear recurrence with
polynomial coefficients. This procedure is called
{\tt SchemeRecurrence}. The reader should look up the on-line help.
 
{\tt SchemeF} does what Scheme does, but more generally,
by trying {\tt Scheme} all on the images of $\P$ under
the dihedral group consisting of inverse, and reverse.
 
This paper's website 
{\tt http://www.math.temple.edu/\Tilde zeilberg/WILF.html}
has numerous sample input and output files.
 
{\bf The Maple package } {\tt HERB}
 
This is the non-rigorous counterpart of {\tt WILF}.
There a scheme is found empirically by checking
$(Julian)$ for small $n$, and extrapolating.
It was written before the rigorous {\tt WILF}, and
helped a lot in developing the latter. 
 
{\bf Future Directions}
 
Prefix Schemes are equivalent to Suffix Schemes, but it should
be possible to have {\it mixed schemes}. Also one may be even
more creative in partitioning $A_\sigma$. Hopefully,
a more general notion of {\it scheme} would be more
successful. Also, it should be possible to
{\it empirically} guess {\it generating functions}
(or perhaps, {\it redundant generating functions} (in the
sense of MacMahon)), for the $A_\sigma$, 
in $Redu \cup Expa$, which, once guessed, are rigorously
provable. This, in particular, would entail a
{\it constant-term expression} for $A(n; \P)$, which,
by using the WZ method, should lead to a {\it rigorous}
derivation of the recurrence guessed by procedure
{\tt SchemeRecurrence}  in {\tt WILF}.
 
{\bf REFERENCES}
 
[K] Donald E. Knuth, ``{\it The Art of Computer Programming}'',
vol.3,  Addison-Wesley, Reading, MA, 1973. Second Edition 1997.
 
[L] Levi Ben Gerson, ``{\it Sefer Ma'asei Khosev}'', 1321, Orange, France.
[Extant printed version in: Gerson Lange, Frankfurt: Louis Golde, 1909.]
 
[We] Julian West, {\it Generating trees and Catalan and Schr\"oder numbers},
Discrete Mathematics {\bf 146}(1995), 247-262.
 
[Wi] Herbert S. Wilf, {\it What is an Answer?}, Amer. Math. Monthly
{\bf 89}(1982), 289-292.
 
\bye